\documentclass[12pt]{article}
\usepackage{amssymb}
\usepackage{epsfig}
\usepackage{mathrsfs}

\makeatletter
\ifnum\@ptsize=0 \advance\hoffset by -0.3cm \fi
\ifnum\@ptsize=2 \advance\hoffset by 0.5cm \fi

\newtheorem{satz}{Satz}[section]

\newtheorem{corollary}[satz]{Corollary}
\newtheorem{definition}[satz]{Definition}
\newtheorem{example}[satz]{Example}

\newtheorem{lemma}[satz]{Lemma}
\newtheorem{proposition}[satz]{Proposition}
\newtheorem{remark}[satz]{Remark}

\newtheorem{theorem}[satz]{Theorem}
\newtheorem{thevarthm}[satz]{\varthmname}

\newenvironment{varthm*}[1]{\trivlist\item[]{\bf #1.}\it}
   {\endtrivlist}

\newcommand\beginproof[1]{%
   \trivlist\item[\hskip\labelsep{\em #1.}]}
\newcommand\proof{\beginproof{Proof}}
\newcommand\proofof[1]{\beginproof{Proof of #1}}
\def\endproof{\hspace*{\fill}\endproofsymbol\endtrivlist}
\def\endproofsymbol{\frame{\rule[0pt]{0pt}{6pt}\rule[0pt]{6pt}{0pt}}}

\let\oldlarge=\large
\let\large=\centering

\renewcommand\@seccntformat[1]{\csname the#1\endcsname.\enspace}

\renewcommand\emptyset{\varnothing}  
\renewcommand\ge{\geqslant}  
\renewcommand\le{\leqslant}  
\renewcommand\epsilon{\varepsilon}

\renewcommand\phi{\varphi}

\renewcommand\({\left(}
\renewcommand\){\right)}

\newcommand\calo{{\mathcal O}}
\newcommand\cale{{\mathcal E}}

\newenvironment{items}
   {\list{\labelitemi}{
      \parsep=0cm \itemsep=0cm \topsep=0cm \partopsep=0.5\baselineskip
      \def\makelabel##1{\hss\llap{\rm##1}}}}
   {\endlist}

\newcommand\norm[1]{\left\|#1\right\|}
\newcommand\twoconditions[2]{_{\scriptstyle#1\atop\scriptstyle#2}}
\newcommand\restr[1]{\big|_{#1}}
\newcommand\rounddown[1]{\left\lfloor#1\right\rfloor}
\newcommand\set[1]{\left\{\,#1\,\right\}}
\newcommand\with{\ \vrule\ }
\newcommand\lra{\longrightarrow}

\newcommand\be{\begin{eqnarray*}}
\newcommand\ee{\end{eqnarray*}}
\newcommand\eqnref[1]{(\ref{#1})}

\newcommand\eps{\epsilon}
\newcommand\vect[1]{\left(\begin{array}{c}#1\end{array}\right)}

\newcommand\Q{\mathbb Q}
\newcommand\bbQ{\mathbb Q}
\newcommand\bbP{\mathbb P}
\newcommand\R{\mathbb R}

\newcommand\bbZ{\mathbb Z}
\newcommand\adj{_{\rm adj}}
\newcommand\Sinv{S^{-1}}

\newcommand\newop[2]{\def#1{\mathop{\rm #2}\nolimits}}
\newop\supp{supp}
\newop\SB{SB}
\newop\SC{SC}
\newop\Bstab{B_{+}}
\newop\End{End}
\newop\Null{Null}
\newop\Neg{Neg}
\newop\Nef{Nef}
\newop\Amp{Amp}
\newop\Face{Face}
\newop\BigCone{Big}
\newop\inter{int}
\newop\NS{NS}
\newop\voll{vol}
\newop\Pic{Pic}

\newcommand\vecspan[1]{\left\langle#1\right\rangle}
\newcommand\orth{^\perp}
\newcommand\nonneg{^{\geqslant 0}}
\newcommand\pos{^{>0}}
\newcommand\relint{\mbox{\rm rel.int.}}
\newcommand\I{\mathcal{I}}
\newcommand\NR{N^1_\R}

\newcommand{\zd}{Zariski decomposition\ }

\newcommand{\ep}{\ensuremath{\epsilon}}
\newcommand{\sr}{\ensuremath{\sqrt{45+78{\epsilon}+49{\epsilon}^2}}}

\newcommand{\hh}[3]{\ensuremath{h^{#1}\left(#2,#3\right)}}
\newcommand{\HH}[3]{\ensuremath{H^{#1}\left(#2,#3\right)}}
\newcommand{\vl}[1]{\ensuremath{{\rm vol}\left(#1\right)}}
\newcommand{\vol}[2]{\ensuremath{{\rm vol}_{#1}\left( #2 \right) } }
\newcommand{\si}[2]{\ensuremath{ \left( { \left( #1 \right) }^{#2}\right) } }

\newcommand{\zj}[1]{\ensuremath{ \left( #1 \right) }}

\newcommand{\NN}{\ensuremath{\mathbb N}}
\newcommand{\PP}{\ensuremath{\mathbb P}}
\newcommand{\QQ}{\ensuremath{\mathbb Q}}
\newcommand{\RR}{\ensuremath{\mathbb R}}

\newcommand{\OO}{\ensuremath{\mathcal O}}
\newcommand{\II}{\ensuremath{\mathcal I}}

\begin{document}

\title{\oldlarge\bfseries Zariski chambers, volumes, and stable base loci}
\author{\normalsize Th.~Bauer, A.~K\"uronya, T.~Szemberg}
\date{\normalsize Version of \today}
\maketitle

\thispagestyle{empty}


\section*{Introduction}

   In this paper we consider certain asymptotic invariants of
   linear systems on algebraic surfaces.
   Originally we were interested in understanding how
   the volume of a line bundle and its stable base locus
   behaves with respect to (small) perturbations of the considered
   bundle. These are asymptotic variants of questions studied
   classically in algebraic geometry. The asymptotic approach
   emerged only recently and quickly gained considerable
   interest.

   Let $X$ be a smooth projective variety of dimension $n$, $L$ 
   a line bundle on $X$. The Riemann--Roch problem is concerned 
   with the study of the behaviour of $\hh{0}{X}{kL}$ as a function
   of $k$. The exact determination of $h^0$ of
   a line bundle is difficult in general. The groups in question 
   typically grow like $k^{\dim X}$ so we can introduce their 
   asymptotic counterpart, the notion of the {\it volume}
   of a line bundle $L$
   $$
      {\rm vol}_X(L)\stackrel{\rm
      def}{=}\limsup_k{\frac{\hh{0}{X}{L^{\otimes k}}}{k^n/n!}}\ ,
      \ 
   $$
   which behaves much more nicely  in many important cases. The volume was first
   introduced by Cutkosky (notably in his proof of non-existence of
   Zariski decomposition for divisors on higher dimensional
   varieties) and subsequently studied by Demailly, Ein, Lazarsfeld and
   others. 

   The concept readily extends --- via homogeneity --- to $\QQ$-divisors and 
   quite generally, it enjoys many useful properties. It has been established 
   recently (in \cite{PAG}) that on an irreducible projective variety of 
   dimension $n$ the volume defines a continuous function on the 
   N\'eron--Severi space.  Also, the volume  is log-concave and 
   homogeneous of degree $n$.


   Still, there have been few instances worked out in the
   literature so far. Here we show that the volume function on the
   cone of big divisors on an algebraic surface is piecewise polynomial
   (more precisely, the big cone splits into subcones on which the volume
   function is quadratic).



   For the base loci of linear series the situation is similar.
   Whereas  determining the
   base locus of a given line bundle is quite difficult in
   general, its asymptotic version is somehow
   easier to study and exhibits better behavior. Given a line
   bundle $L$, the stable base locus $\SB(L)$ of $L$ is
   the intersection of the base loci of the linear series
   $|kL|$ for all positive integers $k$.
   More generally, we consider the stable base loci of
   $\Q$-line bundles $L$ by passing to an integral multiple
   of $L$; this is well-defined, since the stable base locus
   is invariant under taking
   multiples (i.e.\ tensor powers) of a line bundle.
   These objects were recently studied by Nakamaye (\cite{Nak1},
   \cite{Nak2}). He showed in particular that stable base loci of
   slightly perturbed nef divisors remain constant in a small
   neighborhood of the perturbed divisor. Here we investigate in
   more detail the regions (in the big cone) where
   the stable base loci remain constant. These regions in the case
   of surfaces turn out to be
   convex polyhedral
   subcones.

   Both problems thus lead to a partition of the big cone
   into suitable subcones, and it is natural to ask whether the
   partitions agree. Somewhat surprisingly, we show that this is
   indeed the case and that in fact both
   problems are closely related to the variation of the Zariski
   decomposition, which is an interesting problem quite on its
   own.

   More precisely, knowing the Zariski decomposition of a $\QQ$-divisor 
   provides a quick way to determine both the volume and the stable
   base locus of the divisor. Therefore the description of the regions
   where the support of the negative part of Zariski decompositions
   is constant will settle both questions we originally studied.

   Our main result is the following

%
%

\begin{varthm*}{Theorem}
   Let $X$ be a smooth projective surface. Then there is a
   locally finite decomposition of the big cone of $X$ into
   rational locally polyhedral subcones
   such that the following holds:
\begin{items}
\item[(i)]
   In each subcone the
   support of the negative part of the Zariski decomposition of
   the divisors in the subcone is constant.
\item[(ii)]
   On each of the subcones the volume function is given by a
   single polynomial of degree two.
\item[(iii)]
   In the interior of each of the subcones the stable
   base loci are constant.
\end{items}
\end{varthm*}

   In addition,  we work out in detail the volume function
   on  del Pezzo and $K3$ surfaces and explore the
   connections with the action of the Weyl group. In contrast to the case 
   of surfaces we establish the following theorem.

   \begin{varthm*}{Theorem} 
     For every $n\geq 3$ there exists a smooth projective 
     variety of dimension $n$ such that its corresponding volume 
     function is not locally polynomial.
   \end{varthm*}

   The organization of the paper goes as follows. In Section 1 we consider the 
   problem of variation of Zariski decompositions in the big cone of a smooth 
   projective surface. We establish the part of the Theorem which regards
   Zariski decompositions. Along the way we prove that in the interior of the 
   big cone the nef cone is actually locally polyhedral which gives a strengthening
   of a result of Campana and Peternell (\cite{CamPet90}) on the geometry of 
   the nef boundary. 
 
   Section 2 deals with the applications of the results in Section 1 to the
   description of stable base loci and destabilizing numbers. We prove that 
   on surfaces all destabilizing numbers --- with the possible exception of the 
   largest one --- of a big divisor $L$ with respect to an ample divisor $A$ are rational. 
   We give a counterexample to this statement in dimension three.

   We then move  on to describe the volume function on surfaces in Section 3, with
   detailed computations in the case of del Pezzo surfaces. In addition, we 
   investigate the relation between the volume and the Weyl action on K3 surfaces.
   Lastly, we provide an example of smooth $n$-folds ($n\geq 3$) where the volume function is 
   not locally polynomial.

   Finally, Section 4 contains a few somewhat technical  results that are used in the main text.
   We include them for the sake of completeness  as we were not able to find a 
   refererence for them.

   \textbf{Acknowledgements.} Thomas Bauer was supported by the DFG grant 
   BA 1559/4. Tomasz Szemberg was partially supported by the  KBN grant 2P03A 022 17.
   The authors would like to thank Mihnea Popa,
   Zolt\'an Szab\'o and \'Arp\'ad T\'oth for interesting and helpful discussions. 
   Special thanks are due to Igor Dolgachev for  sharing with us his insight   
   on $K3$ surfaces and Lawrence Ein and Mircea Mustata for pointing out mistakes 
   in one of the predecessors of this manuscript (due to the second author).
   Finally, the authors would like to express their gratitude towards Rob
   Lazarsfeld, the thesis advisor of the second author,  for suggesting the 
   cooperation, his support  and  many useful comments.


\section{Zariski decompositions}

   In this section we will prove the following theorem on the
   variation of the Zariski decomposition in the big cone.

\begin{theorem}\label{Zariski chamber theorem}
   Let $X$ be a smooth projective surface. Then there is a
   locally finite decomposition of the big cone of $X$ into
   rational locally polyhedral subcones
   such that in each subcone the
   support of the negative part of the Zariski decomposition of
   the divisors in the subcone is constant.
\end{theorem}

   We will use the following notation. If $D$ is an $\R$-divisor,
   we will write
   $$
      D=P_D+N_D
   $$
   for its Zariski decomposition, and we let
   $$
      \Null(D)=\set{C\with C \mbox{ irreducible curve with } D\cdot C=0}
   $$
   and
   $$
      \Neg(D)=\set{C\with C \mbox{ irreducible component of }
      N_D}\ .
   $$
   Of course $\Neg(D)\subset\Null(P_D)$.

\begin{proposition}\label{prop1}
   Let $X$ be a smooth projective surface and $P$ a big and nef
   $\R$-divisor on $X$. Then there is a
   neighborhood $U$ of $P$ in $\NR(X)$ such that for all divisors
   $D\in U$ one has
   $$
      \Null(D)\subset\Null(P) \ .
   $$
\end{proposition}

\proof
   As the big cone is open, we may
   choose big (and effective) $\R$-divisors $D_1,\dots,D_r$ such that
   $P$ lies in the interior of the cone
   $\sum_{i=1}^r\R^+D_i$.
   We can have $D_i\cdot C<0$ only for finitely many curves
   $C$. Therefore,
   after possibly replacing $D_i$ with $\eta D_i$ for some small
   $\eta>0$, we can
   assume that
   $$
      (P+D_i)\cdot C>0    \eqno (*)
   $$
   for all curves $C$ with $P\cdot C>0$.
   We conclude then from $(*)$ that
   $$
      \Null\(\sum_{i=1}^r \alpha_i (P+D_i)\)\subset\Null(P)
   $$
   for any $\alpha_i>0$. So the cone
   $$
      U=\sum_{i=1}^r \R^+(P+D_i)
   $$
   is a neighborhood of $P$ with the desired property.
\endproof

   Denote by $\I(X)$ the set of all irreducible curves on $X$ with
   negative self-intersection.
   Note that if $D$ is a big divisor, then by the
   Hodge index theorem
   we have $\Null(D)\subset\I(X)$.
   For $C\in\I(X)$ denote by $C\nonneg$
   the half-space
   $\set{D\in\NR(X)\with D\cdot C\ge 0}$
   and by $C\orth$ the hyperplane
   $\set{D\in\NR(X)\with D\cdot C=0}$.

\begin{corollary}\label{cor big nef polyhedral}
   The intersection of the nef cone and the big cone is locally
   polyhedral, i.e., for every $\R$-divisor
   $P\in\Nef(X)\cap\BigCone(X)$ there exists a neighborhood
   $U$ and curves $C_1,\dots,C_k\in\I(X)$ such that
   $$
      U\cap\Nef(X)=U\cap\(C_1\nonneg\cap\dots\cap C_k\nonneg\)
   $$
\end{corollary}

\proof
   Let $U$ be a neighborhood of $P$ as in the proposition.
   We have
   $$
      \BigCone(X)\cap\Nef(X)=\BigCone(X)\cap\bigcap_{C\in\I(X)}C\nonneg
   $$
   and therefore
   $$
      U\cap\Nef(X)=U\cap\bigcap_{C\in\I(X)}C\nonneg
      \eqno(*)
   $$
   For every $C\in\I(X)$ we have either $U\subset C\nonneg$, in
   which case we may safely omit $C\nonneg$ from the intersection in
   $(*)$, or else $U\cap C\orth\ne\emptyset$.
   But by choice of $U$,
   the second option
   can only happen for finitely many curves $C$.
   In fact, $U\cap\Nef(X)=U\cap\bigcap_{C\in\Null(P)}C\nonneg$.
\endproof

   Let $P$ be a big and nef $\R$-divisor on $X$.
   The \textit{face} of $P$ is given by
   $$
      \Face(P)=\bigcap\twoconditions{C\in\I(X)}{P\in C\orth}
         C\orth \cap\Nef(X)  \\
      =  \bigcap_{C\in\Null(P)} C\orth \cap\Nef(X)
   $$
   so that
   \begin{equation}\label{face formula}
      \Face(P)=\Null(P)\orth\cap\Nef(X)
   \end{equation}

   Given a big and nef $\R$-divisor $P$, consider the
   set
   $$
      \Sigma_P=\set{D\in\BigCone(X)\with\Neg(D)=\Null(P)}
      \ .
   $$
   One checks that
   $\Sigma_P$ is a convex cone. It will in general
   be neither open nor closed.
   In Example \ref{two_points} the chamber $\Sigma_{Q_2}$ contains
   the wall spanned by $0$, $L$ and $E_2$ but it doesn't
   contain the ray through $L$ nor the wall spanned by $0$, $L$ and
   $L-E_1$. On the other hand $\Sigma_L$ is open.
   If $A$ is an ample divisor, then
   $\Sigma_A=\Nef(X)$ is closed. This shows that
   all possibilities for the boundary points can happen.

   Our aim is now to show that the cones $\Sigma_P$ provide the
   decomposition that is claimed in the theorem.
   We start with the
   following properties of these cones.

\begin{lemma}\label{lemma chambers}
   Let $P$ and $P'$ be big and nef divisors on $X$.
\begin{items}
\item[(i)]
   $\Sigma_P=\Sigma_{P'}$ if and only if $\Face(P)=\Face(P')$.
\item[(ii)]
   $\Sigma_P\cap\Sigma_{P'}=\emptyset$, if $\Face(P)\ne\Face(P')$.
\item[(iii)]
   $\BigCone(X)$ is the union of the sets $\Sigma_P$.
\end{items}
\end{lemma}

\proof
   (i) and (ii) follow from \eqnref{face formula} plus the fact that
   $\Null(P)\orth\cap\Nef(X)=\Null(P')\orth\cap\Nef(X)$
   implies $\Null(P)=\Null(P')$.

   For (iii),
   given a big $\R$-divisor $D$, we need to show that there is a
   big and nef divisor $P$ such that $\Neg(D)=\Null(P)$.
   Let $\Neg(D)=\set{C_1,\dots,C_k}$, and take any ample divisor
   $A$.
   We claim that a divisor $P$ as required can be constructed
   explicitly in the form
   $A+\sum_{i=1}^k\lambda_i C_i$ with suitable non-negative
   rational numbers
   $\lambda_i$. In fact, the conditions to be fulfilled are
   $$
      \(A+\sum_{i=1}^k\lambda_i C_i\)\cdot C_j=0
        \quad\mbox{ for } j=1,\dots,k \ .
   $$
   This is a system of linear equations with negative definite
   coefficient matrix
   $(C_i\cdot C_j)$, and Lemma \ref{lemma matrix} guarantees
   that all components $\lambda_i$
   of its solution are non-negative.
   In fact all $\lambda_i$'s  must be positive as $A$ is ample.
\endproof

   The following proposition gives a useful characterization
   of the loci where two or more faces meet.

\begin{proposition}\label{boundary}
   Let $D$ be a big divisor on $X$. Then $D$ is in the boundary of
   some $\Sigma_P$ if and only if $\Neg(D)\neq\Null(P_D)$.
\end{proposition}

\proof
   Let $D\in\partial\Sigma_P$ and let $N_D=\sum_{i=1}^m a_iN_i$. We fix a
   norm $\norm\cdot$ on $\NR(X)$. Then for every small $\eps>0$
   there exists an element $\alpha\in\NR(X)$ with
   $\norm\alpha<\eps$ such that $\Neg(D+\alpha)\neq \Neg(D)$.

   Let $\alpha=\alpha'+\alpha''$ be the
   decomposition induced by the
   direct sum
   $$
     \NR(X)=\vecspan{P_D,N_1,\dots,N_m}\oplus\vecspan{P_D,N_1,\dots,N_m}\orth
   $$
   and observe that $\Neg(D+\alpha')=\Neg(D)$ for $\alpha'$ of small norm.
   So we may assume that $\alpha\in\vecspan{P_D,N_1,\dots,N_m}\orth$.

   Now we take a convex neighborhood $U$ of $P_D$ satisfying
   the condition of Proposition \ref{prop1}.
   Rescaling if necessary we may assume that $P_D+\alpha\in U$.
   Moreover $P_D+\alpha$ is not nef as otherwise we would have
   the Zariski decomposition $D+\alpha=(P_D+\alpha)+N_D$. Hence there
   is a curve $C_\alpha\in \I\setminus\{N_1,\dots,N_m\}$ with
   $(P_D+\alpha)\cdot C_\alpha<0$.
   This implies that for some $t\in[0,1]$ we have
   $(P_D+t\alpha)\cdot C_\alpha=0$.
   So we have $C_\alpha\in\Null(P_D+t\alpha)$, which
   by the choice of $U$ implies $C_\alpha\in\Null(P_D)$.
   This shows the strict inclusion of $\Neg(D)=\{N_1,\dots,N_m\}$ in $\Null(P_D)$.

   For the other direction assume that there is an irreducible
   negative curve $C\in\Null(P_D)\setminus\Neg(D)$.
   One checks that for $\eps>0$ the divisor $D+\eps C$ has
   the Zariski decomposition
   $$
     D+\eps C=P_D+(N_D+\eps C) \ .
   $$
   Letting $\eps$ converge to $0$ this means
   that $D$ can be approximated by divisors from a different
   chamber, so it
   must be in the boundary. (Note that e.g.\
   $D\in\partial\Sigma_{P_D}$.)
\endproof

   We give now a description of the interior of the
   chambers.

\begin{proposition}
   The interior of $\Sigma_P$ is given by
   $$
       \set{ D\in\BigCone(X)\with \Neg(D)=\Null(P)=\Null(P_D) }
   $$
\end{proposition}

\proof
   As the chambers $\Sigma_{P}$ are disjoint, the interior of
   $\Sigma_P$ consists of the points that are not on the boundary of
   some chamber.
\endproof

We turn now to the description of the closure of the chambers $\Sigma_P$.

\begin{corollary}\label{closure}
   Let $P$ be a nef and big divisor. Then
   $$
     \BigCone(X)\cap\overline{\Sigma}_P=
       \set{D\in\BigCone(X) \with
       \Neg(D)\subset\Null(P)\subset\Null(P_D)}
   $$
\end{corollary}

\proof
   Let $D$ be in the set on the right. Then $P_D\in\Face(P)$ because
   $$
    P_D\in \Null(P_D)\orth\cap\Nef(X) \subset
    \Null(P)\orth\cap\Nef(X)=\Face(P).
   $$
   Hence there is a sequence of divisors $Q_n$ in the relative interior of $\Face(P)$
   converging to $P_D$. In particular for every $n$ we have $\Null(Q_n)=\Null(P)$.

   On the other hand, as $\Neg(D)\subset\Null(P)$ there exists a sequence $N_n$
   of negative definite divisors converging to $N_D$ such that
   $\Neg(N_n)=\Null(P)$.
   (Just add to $N_D$ small fractions of curves $C\in\Null(P)\setminus\Neg(D)$ and
   apply as usual Lemma \ref{missing lemma}.)

   Putting things together, we obtain a sequence $D_n=Q_n+N_n$ of divisors
   converging to $D$. As $D_n=Q_n+N_n$ is by construction the Zariski decomposition,
   the divisors $D_n$ lie all in the interior of $\Sigma_P$.
   This implies that $D\in\overline{\Sigma}_P$.

   For the other direction, assume that
   $D\in\BigCone(X)\cap\overline{\Sigma}_P$
   and let $D_n$ be a sequence of divisors in the interior of
   $\Sigma_P$
   converging to $D$. This implies that for every $n$
   $$
      \Neg(D_n)=\Null(P)
      \mbox{ and }
      \Null(P_{D_n})=\Null(P).
   $$
   In other words
   $$
      P_{D_n}\in \Null(P)\orth
      \quad\mbox{and}\quad
      N_{D_n}\in \vecspan{\Null(P)}
   $$
   for every $n$. As both spaces are closed and orthogonal
   to each other, we obtain
   $$
      P_{D_n}\longrightarrow P_D \in \Null(P)\orth
      \quad\mbox{and}\quad
      N_{D_n}\longrightarrow N_D \in \vecspan{\Null(P)}
   $$
   Dualizing the first condition and using Lemma
   \ref{curves in negative}
   for the second we arrive at
   $$
      \Null(P)\subset \Null(P_D)
      \quad\mbox{and}\quad
      \Neg(D)\subset \Null(P).
   $$
\endproof

   We give now an explicit description of the chambers
   $\Sigma_P$.
   Here we use the notation
   $V\nonneg(M)$ for the subcone of $\NR(X)$ generated
   by a subset $M\subset\NR(X)$
   and $V\pos(M)$ for its interior.

\begin{proposition}\label{prop chamber closures}
   $$
   \BigCone(X)\cap\overline\Sigma_P=\(\BigCone(X)\cap\Face(P)\)+V\nonneg\(\Null(P)\)
   $$
\end{proposition}

\proof
   We will use the characterization of the closure of $\Sigma_P$
   we have just obtained. Take a big divisor $D$ for which
   \[
      \Neg(D)\subset \Null(P)\subset \Null(P_D)
   \]
   holds. Then $P_D\in \Null(P_D)\orth \subset \Null(P)\orth$, hence
   $P_D\in\Face(P)\cap\BigCone(X)$. On the other hand
   $\Neg(D)\subset \Null(P)$ implies $N_D\in
   V\nonneg\(\Null(P)\)$, so that we have $D\in
   \(\BigCone(X)\cap\Face(P)\)+V\nonneg\(\Null(P)\)$.

   Going in the
   other direction, pick a big divisor
   $D\in\(\BigCone(X)\cap\Face(P)\)+V\nonneg\(\Null(P)\)$. This
   latter is an orthogonal decomposition. Let
   $Q$ and $M$ be
   the components of $D$ in its respective parts. Then $Q\cdot M=0$,
   $Q$ is nef, and $M$ is a negative definite divisor. Therefore $D=Q+M$ is the
   Zariski decomposition of $D$ (by the uniqueness of Zariski
   decompositions). Evidently, $\Neg(D)=\Neg(M)\subset \Null(P)$
   (the inclusion coming from Lemma \ref{curves in negative}),
   and $Q=P_D\in\Face(P)$, which
   implies $\Null(P)\subset\Null(P_D)$.
\endproof

\begin{proposition}
   The interior of the chamber $\Sigma_P$ is equal to
   \[
      \relint\Face(P)+V\pos(\Null(P)) \ .
   \]
\end{proposition}

\proof
   Let $D$ be a big divisor in the interior of $\Sigma_P$. Then
   \[
      \Neg(D)=\Null(P)=\Null(P_D) \ ,
   \]
   hence $P_D\in \relint\Face(P)$ and $N_D\in V\pos(\Null(P))$.

   On the other hand, if $D\in\relint\Face(P)+V\pos\Null(P)$ has the decomposition
   \[
      D=Q+M
   \]
   with respect to the linear subspaces generated by
   $\relint\Face(P)$ and $V\pos(\Null(P))$, then $M$ is a
   negative definite
   divisor (or zero), $Q\cdot M=0$, and $Q$ is nef, therefore $D=Q+M$ is
   again the Zariski decomposition of $D$. As $Q\in\relint\Face(P)$,
   we have $\Null(P)=\Null(Q)=\Null(P_D)$. Also, $\Neg(D)=\Neg(M)=\Null(P)$.
\endproof

   It follows in particular that the chambers are locally polyhedral.
   Within the big cone, the situation is even better:

\begin{proposition}
   If $\Face(P)$ is contained in $\BigCone(X)$,
   then $\Sigma_P$ is polyhedral.
\end{proposition}

\proof
   Take a hyperplane $H$ in $\NR(X)$ cross-secting $\Nef(X)$.
   Then $\Face(P)\cap H$ is compact.
   Corollary \ref{cor big nef polyhedral} implies that
   every divisor $D\in\Face(P)\cap H$ has
   an open neighborhood $U$ such that $U\cap\Face(P)\cap H$ is
   polyhedral. By compactness we conclude that $\Face(P)\cap H$
   itself is polyhedral, and this implies that $\Face(P)$ is
   polyhedral.
\endproof

   On the other hand one can easily have faces that are not
   polyhedral:

\begin{example}\rm
   Take a surface $X$ with infinitely many $(-1)$-curves
   $C_1,C_2,\dots$, and blow it up at a point that is not
   contained in any of the curves $C_i$. On the blow-up
   consider the exceptional divisor $E$ and the proper transforms
   $C'_i$.
   Since the divisor $E+C'_i$ is negative definite,
   we can proceed as in the proof of Lemma \ref{lemma chambers}
   to construct for every index $i$
   a big and nef divisor $P_i$ with
   $\Null(P_i)=\set{E,C'_i}$, and also a divisor
   $P$ such that $\Null(P)=\set{E}$.
   But then $\Face(P)$ meets contains
   all faces $\Face(P_i)$, and
   therefore it is
   not polyhedral.
\end{example}

   In order to prove local finiteness of the chamber
   decomposition we
   will also make use of
   the following statement about Zariski
   decompositions.

\begin{lemma}\label{lemma add ample}
   If $D$ is a big divisor and $A$ is an ample divisor, then for
   all $\lambda\geq 0$
   $$
      \Neg(D+\lambda A)\subset\Neg(D) \ .
   $$
\end{lemma}

\proof
   The idea is to proceed by a (finite) induction on
   the number of elements in $\Neg(D)$.

   If $\Neg(D)=\emptyset$, then $D=P_D$ is nef, hence
   $D+A$ is ample and $\Neg(D+A)$ is an empty set as well.

   Assume now for $r\geq 1$ that the Lemma holds for all big divisors
   with at most $r-1$ irreducible components in the negative part
   of their Zariski decompositions and let $D=P_D+N_D$
   be a big divisor with $N_D=\sum_{i=1}^r a_iN_i$. We will prove
   the following
\begin{quote}
   {\bf Claim.}
   There exists a positive number $\eps_0>0$ and affine-linear functions
   $f_1,\dots,f_r:\R\to\R$
   such that for $\eps$ with $0\le\eps\le\eps_0$
   the Zariski decomposition of $L+\eps A$ is
   $$
      \(P+\eps A+\sum_{i=1}^r(a_i-f_i(\eps))N_i\)
      +\sum_{i=1}^r f_i(\eps)N_i
   $$
   (the expression in brackets being the positive part)
   and such that $\eps_0$ is a zero of one of the functions $f_i$.
\end{quote}

   We show first that the above claim suffices in order to
   complete the proof of the Lemma. Indeed,
   for $\lambda\in [0,\eps_0]$ the statement follows from the
   claim, whereas for $\lambda\geq\eps_0$ we have
   $$
      D+\lambda A=(D+\eps_0 A)+(\lambda-\eps_0)A
   $$
   and the induction hypothesis applies to $D+\eps_0 A$.

   Turning to the proof of the claim
   consider the $\R$-divisor
   $$
      P'=P+\eps A+\sum_{i=1}^r(a_i-x_i)N_i \ .
   $$
   The divisor $L+\eps A$ has the Zariski decomposition
   $$
      P'+\sum_{i=1}^r x_iN_i
   $$
   if the following conditions are satisfied
   \begin{eqnarray*}
      & (1) & 0\le x_i\le a_i \quad\mbox{ for all i} \\
      & (2) & P'\cdot N_i=0 \\
      & (3) & P' \mbox{ is nef}
   \end{eqnarray*}
   Note that (3) follows from (1) and (2), since $P+\eps A$ is ample.
   Condition (2) is equivalent to the linear system of equations
   $$
      S\cdot\vect{x_1\\ \vdots \\ x_r}=\eps\vect{AN_1\\ \vdots\\ AN_r}
         + \vect{NN_1\\ \vdots \\ NN_r}
   $$
   where $S$ denotes the intersection matrix of $N$. Since $S$ is negative definite,
   the system has the unique solution
   $$
      \vect{x_1\\ \vdots \\ x_r}=\eps\Sinv\vect{AN_1\\ \vdots\\ AN_r}
         +\vect{a_1\\ \vdots\\ a_r}
   $$
   So the $x_i$ are linear functions $f_i$ of the parameter $\eps$.
   Furthermore, we find $x_i\le a_i$, because by Lemma \ref{lemma matrix}
   below all the entries of the matrix $\Sinv$ are $\le 0$.
   Thus the conditions (1) and (2) will be satisfied if we choose
   $\eps_o$ as the smallest zero of the functions $f_i$.
\endproof

\begin{proposition}\label{prop decomposition locally finite}
   The decomposition of the big cone into the chambers $\Sigma_P$
   is locally finite.
\end{proposition}

\proof
   Denote by $\Amp(X)$ the ample cone of $X$.
   Every big divisor has an open neighborhood
   in $\BigCone(X)$ of the form
   $$
      D+\Amp(X)
   $$
   for some big divisor $D$.
   (In fact, given a big divisor $D_0$, there is an ample divisor
   $A$ such that $D_0-A$ is still big, and hence
   $D_0\in(D_0-A)+\Amp(X)$.)
   Lemma \ref{lemma add ample} implies that only finitely many
   chambers $\Sigma_P$ can meet this neighborhood.
\endproof

   We give now the proof of the theorem.

\proofof{Theorem \ref{Zariski chamber theorem}}
   According to Lemma~\ref{lemma chambers} the subcones $\Sigma_P$ yield a
   decomposition of the big cone. By definition the support of
   the negative part of the Zariski decomposition is constant in
   $\Sigma_P$.
   Prop.~\ref{prop chamber closures} implies that the $\Sigma_P$
   are locally polyhedral, and
   Prop.~\ref{prop decomposition locally finite}
   completes the proof.
\endproof

   As a first application, we
   now show that Zariski decompositions in the big cone are
   continuos.
   The proof uses the local finiteness of the chamber
   structure in an essential way.

\begin{proposition}
   Let $(D_n)$ be a sequence of big divisors converging in
   $\NR(X)$ to a big divisor $D$.
   If $D_n=P_n+N_n$ is the Zariski decomposition of $D_n$, and
   if $D=P+N$ is the Zariski decomposition of $D$, then
   the sequences $(P_n)$ and $(N_n)$ converge to $P$ and $N$
   respectively.
\end{proposition}

\proof
   We consider first the case where all $D_n$ lie in a fixed
   chamber $\Sigma_P$. In that case we have by definition
   $\Neg(D_n)=\Null(P)$
   for all $n$, so that
   $$
      N_n\in\vecspan{\Null(P)}
   $$
   and hence $P_n\in\Null(P)\orth$.
   As
   $$
      \NR(X)=\Null(P)\orth\oplus\vecspan{\Null(P)}
   $$
   we find that both sequences $(P_n)$ and $(N_n)$ are
   convergent.
   The limit class $\lim P_n$ is certainly nef.
   Let $E_1,\dots,E_m$ be the curves in $\Null(P)$. Then
   every $N_n$ is of the form $\sum_{i=1}^m a_i^{(n)}E_i$ with
   $a_i^{(n)}>0$. Since the $E_i$ are
   numerically independent, it follows that $\lim N_n$ is
   of the form $\sum_{i=1}^m a_iE_i$ with $a_i\ge 0$, and hence
   is either negative definite or zero. Therefore
   $D=\lim P_n+\lim N_n$ is actually the Zariski decomposition of
   $D$, and by uniqueness the claim is proved.

   Consider now the general case where the $D_n$ might lie in
   various chambers. Since the decomposition into chambers is
   locally finite, there is a neighborhood of $D$ meeting only
   finitely many of them. Thus there are finitely many
   big and nef divisors $P_1,\dots,P_\ell$ such that
   $$
      D_n\in\bigcup_{i=1}^\ell\Sigma_{P_i}
   $$
   for all $n$.
   So we may decompose the sequence $(D_n)$ into finitely many
   subsequences to which the case above applies.
\endproof


\section{Base loci}

   In this section we study stable base loci as considered recently
   by Nakamaye \cite{Nak1}, \cite{Nak2}.

   For an integral divisor $D$ denote by $\SB(D)$ the
   \textit{stable base locus} of  $D$,
   i.e., the intersection of the base loci of the linear series
   $|kD|$ for all positive integers $k$.
   More generally, we will consider the stable base loci of
   $\Q$-divisors $D$ by passing to an integral multiple
   of $D$; this is well-defined, since the stable base locus is invariant under taking
   multiples (i.e.\ tensor powers) of a line bundle.

   In \cite{AIBL} Ein et al.\ introduced a related
   notion of \textit{stabilized base locus} $\Bstab(M)$ defined as
   $$
      \Bstab(M)=\SB(M-A)
   $$
   for (an arbitrary) sufficiently small ample $\Q$-divisor $A$.
   This notion has the advantage of being independent of the
   numerical equivalence class of a divisor, so it can be safely
   studied in the space $\NR(X)$. Moreover,
   since the definition depends on small perturbations of $D$,
   it extends in a natural manner to $\R$-divisors. Clearly, one
   has always inclusions
   $$
      \Bstab(D+A)\subset \Bstab(D)\subset \Bstab(D-A).
   $$
   An $\R$-divisor $D$ is called \textit{stable} if equalities
   hold, i.e.,
   $$
      \Bstab(D+A) = \Bstab(D) = \Bstab(D-A)
   $$
   for all sufficiently small ample $\R$-divisors $A$.
   A $\Q$-divisor $D$ is stable if and only if
   $\SB(D)= \SB(D\pm A)$
   for $A$ ample and sufficiently small.
   All divisors in $\BigCone(X)$ which are not stable are called
   \textit{instable}. (Our definition, although stated differently,
   agrees with definition 1.25 of \cite{AIBL}.)

\begin{definition}\rm
   If $D$ is a stable $\R$-divisor, then its \emph{chamber of
   stability} is defined as the set
   $$
      \SC(D)=\{D'\in\BigCone(X) | \Bstab(D')=\Bstab(D)\}
      \subset\NR(X).
   $$
\end{definition}

   We will show in this section
   that, somewhat surprisingly, the chambers of stability
   agree essentially with
   the Zariski chambers
   of Theorem \ref{Zariski chamber theorem}.

\begin{theorem}\label{Stable chamber theorem}
   Let $D$ be a stable big $\R$-divisor. Then
   $$
      \inter\SC(D)=\inter\Sigma_{P_D}.
   $$
\end{theorem}

   As it suffices to prove the Theorem for $\Q$-divisors
   we assume from now on that $D$ is a $\Q$-divisor
   and we work with the stable base loci.

   Let $C$ be an irreducible curve on $X$ that is contained in
   $\SB(D)$. We will say that $C$ is
   a \textit{bounded base component} of $D$, if there is a
   constant $p$ such that the coefficient of
   $C$ in the base divisor of the linear series $|kD|$ is less
   than $p$ for all integers $k$ such that $kD$ is an integral
   divisor. We will call $C$ an \textit{unbounded base component}
   otherwise. One checks that this notion is invariant under
   taking multiples of $D$.
   By work of
   Cutkosky and Srinivas \cite{CutSri93} one knows that a bounded
   base component in fact appears with periodic coefficients in the base
   divisors of the linear series $|kD|$ for large $k$.

   We show:

\begin{proposition}\label{prop linear}
   Let $D$ be a big $\Q$-divisor
   and $A$ an ample $\Q$-divisor
   on a smooth projective surface $X$.
   If $C$ is a bounded base component of $D$, then it is an
   unbounded base component of $D-\eps A$ for all rational numbers $\eps>0$.
\end{proposition}

\proof
   Let $D=P_D+N_D$ be the Zariski decomposition of $D$. Then $C$ is a base component
   of $P_D$, since the components of $N_D$ are unbounded base components of $D$
   according to \cite{Zar62}, Theorem 8.1.
   Since $P_D$ is nef and $C$ is a stable base component, we know from
   \cite{Zar62}, Theorem 9.1., that one has $P_D\cdot C=0$,
   so that
   $$
      (P_D-\eps A)C<0 \ .
   $$
   Therefore $C$ is an unbounded base component of $P_D-\eps A$.
   Note that \cite{Zar62}, Corollary 7.2, implies that we have
   $$
      h^0(k(D-\eps A))=h^0(k(P_D-\eps A))
   $$
   for all integers $k$ such that the bundles in question are integral,
   since $(P_D-\eps A)N_i<0$ for all components $N_i$ of $N_D$.
   So $C$ is an unbounded base component of $D-\eps A$.
\endproof

   The proposition above
   gives the following corollary.

\begin{corollary}\label{bounded destabilizing}
   If
   $D$ has a bounded base component, then $D$ is instable.
\end{corollary}

\proof
   Let $C$ be a bounded base component of $D$, and
   let $A$ be an ample divisor.
   For rational numbers $\beta> 0$,
   Proposition \ref{prop linear} implies that $C$ is an unbounded
   base component of $D-\beta A$. In particular, $C$ is contained
   in $\SB(D-\beta A)$. Suppose now that $C$ is also contained in
   $\SB(D+\alpha A)$ for some rational number $\alpha> 0$.
   If $C$ is an unbounded base component of $D+\alpha A$, then it is clearly
   an unbounded base component of $D$ as well, which contradicts the hypothesis.
   So $C$ must be a bounded base component of
   $D+\alpha A$. But then it is an unbounded base component of $D$
   by Proposition \ref{prop linear}, a contradiction again.
\endproof

   The following proposition implies then Theorem \ref{Stable
   chamber theorem}.

\begin{proposition}\label{sb neg}
   If $D$ is a stable $\Q$-divisor, then
   $$
      \Neg(D)=\SB(D) \ .
   $$
\end{proposition}


\proof
   For any rational number $\lambda$ one has
   \be
      \lefteqn{\Neg(D-\lambda A)} \\
      & \subset & \SB(D-\lambda A) \\
      & \subset & \Neg(L-\lambda A)\cup\set{\mbox{bounded base components of $D-\lambda A$}}
   \ee
   The last inclusion here follows from the fact that
   the base components of $P_D$ are bounded, because $P_D$ is nef and big
   (see \cite[Theorem 10.1]{Zar62}),
   plus the fact that the stable base locus does not contain any
   isolated points (see \cite[Theorem 6.1]{Zar62}).
   If $\lambda$ is small enough, then
   every stable base component must be unbounded by Corollary~\ref{bounded destabilizing},
   and the claim follows.
\endproof

\begin{remark}\rm
   Note that using the above proposition one can define stable base loci for stable
   $\R$-divisors.
\end{remark}

\begin{remark}\rm
   It is essential to take interiors in the statement of Theorem
   \ref{Stable chamber theorem} as the boundaries of the chambers
   may differ. For example, if $A$ is an ample divisor, then
   $\SC(A)$ is the open ample cone, whereas $\Sigma_A$ is its
   closure -- the nef cone.
\end{remark}

   We conclude this section by showing that in the case of higher
   dimensional varieties, chambers of stability need not to be
   rational. To this end it is convenient to introduce the
   following notion.

\begin{definition}
   \rm
   Let $L$ be a big line bundle and $A$ an ample line bundle
   on a smooth projective variety $X$.
   A positive real number $\lambda$ is a \textit{destabilizing number} of $L$
   relative to $A$, if
   $$
      \SB(L-\alpha A)\subsetneq\SB(L-\beta A)
   $$
   for all rational numbers $\alpha,\beta$ with $\alpha<\lambda<\beta$.
\end{definition}

   It is clear that the inclusion
   $\SB(L-\alpha A)\subset\SB(L-\beta A)$
   holds for any rational numbers $\alpha<\beta$.
   So the destabilizing numbers are by definition those real numbers,
   where the stable base locus strictly increases, i.e., where one passes
   from one chamber of stability to the another one.
   Note that the smallest and the biggest destabilizing number
   of $L$ relative to $A$ can be conveniently
   characterized by the conditions of nefness
   and bigness:

\begin{remark}\label{remark smallest biggest}
   \rm
   (a) If $L$ is ample, then the smallest destabilizing number of $L$ relative
   to $A$ is the number
   $$
      \sigma=\sup\set{\lambda\with L-\lambda A \mbox{ is nef}} \ ,
   $$
   In fact:
   If $\lambda$ is smaller than $\sigma$, then
   $L-\lambda A$ is ample, so that
   $\SB(L-\lambda A)$ is empty;
   on the other hand, if $\lambda$ is bigger than $\sigma$, then
   we find a curve in $\SB(L-\lambda A)$ accounting for the fact
   that the bundle is not nef.

   (b) The biggest destabilizing number of $L$ relative to $A$ is the number
   $$
      \sup\set{\lambda\with L-\lambda A \mbox{ is big}} \ .
   $$
   This is clear from the fact that
   $\SB(L-\lambda A)$ becomes all of $X$ as $\lambda$ passes this number.
\end{remark}

   Theorem \ref{Stable chamber theorem} and the
   rationality statement in Theorem \ref{Zariski chamber theorem}
   imply:

\begin{proposition}
   Let $X$ be a smooth projective surface,
   let $L$ be a big line bundle and $A$ an ample line bundle on $X$.
   Then all destabilizing numbers of $L$ relative to $A$ with the possible
   exception of the biggest one are rational numbers.
\end{proposition}

   We will show now that this need not to be the case for
   higher dimensional varieties. In particular chambers of
   stability need not to be rational on varieties of dimension
   $\geq 3$. To this end we revoke an example
   of Cutkosky \cite[Example 1.6]{cut} that he used to show the nonexistence
   of a Zariski decomposition in higher dimensions. The idea simply is that
   there is an irrational ray in the boundary of the nef cone which is
   still contained in the big cone.

\begin{example}\label{nonrational}\rm
   Let $C$ be an elliptic curve with $\End C=\bbZ$ and let $S=C\times C$.
   Then the nef cone of $S$ is the circular cone
   $$
      \Nef(S)=\left\{ \alpha\in NS(S) | ({\alpha}^2)\geq 0,
      (\alpha\cdot h)\geq 0\right\}
   $$
   where $h$ is any ample class.
   We denote by $\delta$ the diagonal in $S$ and by $f_1$ and $f_2$
   the fibers of the first, respectively second, projection.
   These divisors generate the N\'eron--Severi group
   and their intersection numbers are
   $(f_1.f_2)=(f_1.\delta)=(f_2.\delta)=1$ and
   $f_1^2=f_2^2=\delta^2=0$.

   Let $f=\delta-f_1-f_2$ and let $V=\bbP(\calo_S(f)\oplus\calo_S)$ be the
   projectivized bundle with the natural projection $\pi:V\lra S$.
   We identify $S$ with the zero section of $V$.
   Let $S_i=\pi^{-1}(f_i)$ for $i=1,2$ and let $H$ be an ample divisor on $V$.
   The divisor $D=H+\alpha S_1+\beta S_2$ is ample for arbitrary
   $\alpha,\beta\in\bbQ_{\geq 0}$ since $S_1$ and $S_2$ are nef.

   It follows that for arbitrary $\gamma\in\bbQ_{\geq 0}$ the divisor
   $D(\gamma)=D+\gamma S$
   is big and that its nefness need only to be tested on curves contained in $S$.
   On $S$ we have
   $d(\gamma)=D(\gamma)\restr{S}=H\restr{S}+\alpha f_1+\beta f_2+\gamma f$.

   Since $f$ is not nef, there exists a maximal positive number $\gamma_0$ (depending
   on $\alpha$ and $\beta$) such that $d(\gamma_0)$ is nef but not ample. The
   condition is simply given by $(d(\gamma_0))^2=0$. For general $\alpha$ and $\beta$
   we obtain an irrational value of $\gamma_0$ and from now on we fix such a pair
   of $\alpha$ and $\beta$.

   We showed above that on $V$ the divisor $D(\gamma)$ is nef for $0\leq\gamma\leq\gamma_0$
   and not nef for $\gamma>\gamma_0$. In particular in the plane (in $N^1(V)$) spanned by
   $S$ and $D$ the ray through $D(\gamma_0)$ is irrational.

   Now it is easy to find two ample divisors $L$ and $A$ in this plane such that the
   half-line $L-\lambda A$ meets the given ray at an irrational point.

   To be more specific, let $p$ be such an integer that $p\gamma_0>1$ and
   let $L=pD+\rounddown{p\gamma_0}S$ and $A=D$. Then the first destabilizing number is
   $$ \lambda_1=\frac{p\gamma_0-\rounddown{p\gamma_0}}{\gamma_0}\notin\bbQ.$$
   Indeed, $\lambda_1$ is the value where the half-line $L-\lambda A$ hits our irrational ray.
   Since $L-\lambda_1A$ is big and $\lambda_1$ is irrational it follows from
   the real valued Nakai-Moishezon criterion \cite{CamPet90} that there exists a divisor,
   in this case obviously $S$, such that $(L-\lambda_1A)^2\cdot S=0$. Hence
   $S\subset\SB(L-\lambda A)$ for $\lambda>\lambda_1$ and rational. This shows
   that $\lambda_1$ is a destabilizing number.
\end{example}


\section{Volume of line bundles on surfaces}

   The results on the variation of the Zariski decomposition make it
   possible to describe the behavior of the volume of line bundles on surfaces.
   Let $X$ be a smooth projective surface.
   Recall that the volume of a  line bundle on $X$ is the nonnegative real number
\[
    {\rm vol}_X(L)\stackrel{\rm def}{=}\limsup_k{\frac{\hh{0}{X}{L^{\otimes k}}}{k^2/2}}\ .
\]
   The definition extends immediately to $\QQ$-divisors by the homogeneity of the volume.
   In the case of surfaces the rays in $\NR(X)$ determined by elements of $\II$
   are all extremal rays of the Mori cone of $X$ (\cite{deb}).
   They also determine the nef cone of $X$ in the sense that it is
   enough to test the nefness of an effective $\QQ$-divisor $D$ by checking
   the non-negativity of $D^2$ and of the
   intersection numbers of $D$ with the elements of $\II$.

\begin{proposition}
   Let $D$ be a big integral divisor, $D=P_D+N_D$ the Zariski decomposition of $D$. Then
   \begin{items}
      \item[(i)] $\HH{0}{X}{kD}=\HH{0}{X}{kP_D}$ for all $k\geq 1$ such that $kP_D$ is integral, and
      \item[(ii)] $ \vl{D}=\vl{P_D}=\zj{P_D^2}.$
   \end{items}
\end{proposition}
   By the homogeneity and continuity of the volume we obtain that

\begin{corollary}
   For an arbitrary big $\RR$-divisor $D$ with Zariski decomposition $D=P_D+N_D$ we have
\[
{\rm vol}(D)=\zj{P_D^2}=\zj{D-N_D}^2\ .
\]
   Therefore on a chamber $\Sigma_Q$ on which the support of the negative part is constant
   the volume is given by
   a homogeneous quadratic polynomial.
\end{corollary}

\begin{corollary}
   The volume function
   $
      \voll_X:\NR(X)\ni D\lra\voll(D)\in\R
   $
   of a smooth projective surface is locally polynomial.
\end{corollary}

\subsection{Del Pezzo surfaces}

   We will work out the volume of line bundles on  del Pezzo
   surfaces and describe its connection to the Weyl action on the
   Picard group. As $-K_X$ is ample, del Pezzo surfaces have only
   a finite number of extremal rays. The corresponding set of
   hyperplanes (consisting of divisors perpendicular to them)
   will give a decomposition of the big cone into a finite set of
   polyhedral chambers on each of which we can write down a
   polynomial formula for the volume.

   Let us establish some notation. We denote by $X=Bl_{\Sigma}({\PP}^2)$ the
   blow-up
   of the projective plane at $\Sigma \subseteq {\PP}^2$ where $\Sigma$ consists
   of at most eight points in general position. The exceptional divisors
   corresponding to the points in $\Sigma$ are denoted by $E_1,\dots ,E_r \
   (r\leq 8)$. We denote the pullback of the hyperplane class on ${\PP}^2$ by $L$.
   These divisor classes generate the Picard group of $X$
   and their intersection numbers are:
   $L^2=1$, $(L.E_i)=0$ and $(E_i.E_j)=-{\delta}_{ij}$ for $1\leq i,j\leq r$.
   For each $1\leq r\leq 8$ one  can describe explicitly all extremal rays on $X$ (see \cite{demz}).
   Recall that a class $\alpha\in\NS(X)$ is a {\it root} if
   ${\alpha}^2=-2$ and $\alpha. K_X=0$. We denote the set of roots
   by ${\mathcal R}$.

\begin{proposition}
   With notation as above, the set $\left\{ E^{\perp} | E\in
   {\II} \right\}$ determines the chambers for the volume function.
   More precisely, we obtain the chambers by dividing the big cone
   into finitely many parts by the hyperplanes $E^{\perp}$.
\end{proposition}

   Together with the combinatorial description of $\II$ this gives
   complete information about the volume on $X$.

\proof
   Observe that as the only negative curves on a del Pezzo surface are $(-1)$-curves, 
   the support of every negative divisor consists of pairwise orthogonal curves. 
   This can be seen as follows.  Take a negative divisor $N=\sum_{i=1}^{m}{a_iN_i}$.  
   Then, as the self-intersection matrix of $N$ is negative definite,  for any $1\leq i<j\leq m$ one has
   \[
       0 > (N_i+N_j)^2 = N_i^2+2(N_i\cdot N_j)+N_j^2 = 2(N_i\cdot N_j)-2 \ .
   \]
   As  $N_i\cdot N_j\geq 0$, this can only hold if $N_i\cdot N_j=0$.

   According to Proposition \ref{boundary}, a big divisor $D$ is in the boundary 
   of a Zariski chamber if and only if
   \[
   \Neg(D)\neq \Null(P_D)\ .
   \]
   From Lemma \ref{missing lemma} we see that if $C\in\Null(P_D)-\Neg(D)$ for an 
   irreducible negative curve $C$ then $N_D+C$ forms a negative divisor. 
   By the previous reasoning, this implies that $N_D\cdot C=0$ hence $D\cdot C=0$, 
   that is, $D\in
   C\orth$ as required.

   Going the other way, if $D\in C\orth$ for an irreducible negative curve $C$ 
   then  either $P_D\cdot C=0$, that is, $C\in\Null(P_D)$ or $P_D\cdot C>0$.

   In the first case, $C\not\in\Neg(D)$, as otherwise we would have $N_D\cdot C<0$ 
   and consequently $D\cdot C<0$ contradicting $D\in C\orth$. Therefore $C\in \Null(P_D)-\Neg(D)$ 
   and $D$ is in the boundary of some Zariski chamber.

   In the second case, $D\cdot C=0$ and $P_D\cdot C>0$ imply $N_D\cdot C<0$. 
   From this we see  that $C\in \Neg(D)$ but this would mean $P_D\cdot C=0$ which is again a contradiction.

   The conclusion is that on a surface on which the only negative curves 
   are $(-1)$-curves, a big divisor $D$ is in the boundary of a Zariski chamber 
   if and only if there exists an $(-1)$-curve $C$ with $D\in C\orth$.
\endproof

\begin{example}[Blow-up of two points in the plane]\label{two_points}\rm
   In this case there are  three irreducible negative curves: the two exceptional divisors,
   $E_1$,  $E_2$ and  the pullback of the line through the two blown-up points, $L-E_1-E_2$. 
   As we saw in the previous proposition, the corresponding hyperplanes determine the 
   chamber structure on the big cone.  They divide the big cone into five regions 
   on each of which the support of the negative part of the \zd remains constant.

   In this case the chambers are simply described as the set of divisors 
   that intersect negatively the same set of negative curves. In the picture, 
   $A$ is any ample divisor, $P,Q_1,Q_2$
   are big and nef divisors in the nef boundary (hence necessarily non-ample)
   which are in the relative interiors of the indicated faces. The chambers we obtain 
   in the big cone are $\Sigma_A$ (the nef cone itself), $\Sigma_P$, $\Sigma_{Q_1}$, 
   $\Sigma_{Q_2}$, $\Sigma_P$ and $\Sigma_L$. Observe that apart from the nef cone, 
   the chambers do not contain the nef divisors they are associated to.

   Let $D=aL-b_1E_i-b_2E_2$ be a big $\RR$-divisor. Then one can express 
   the volume of $D$ in terms of the coordinates $a,b_1,b_2$ as follows:
   \[
       \vl{D}= \left\{ \begin{array}{ll} D^2=a^2-b_1^2-b_2^2 & \textrm{ if $D$ is nef, i.e. $D\in\Sigma_A$} \\
       a^2-b_2^2 & \textrm{ if\ } D\cdot E_1<0 \mbox{ and } D\cdot E_2\geq 0 \textrm{ i.e. $D\in \Sigma_{Q_1}$}\\
       a^2-b_1^2 & \textrm{ if\ } D\cdot E_2<0 \mbox{ and } D\cdot E_1\geq 0 \textrm{ i.e. $D\in\Sigma_{Q_2}$}\\
       a^2 & \textrm{ if } D\cdot E_1<0 \mbox{ and } D\cdot E_2<0 \textrm{ i.e. $D\in\Sigma_L$}\\
       2a^2-2ab_1-2ab_2+2b_1b_2 & \textrm{ if\ } D\cdot(L-E_1-E_2)<0\
       \textrm{ i.e. $D\in\Sigma_P$}.
       \end{array} \right.
    \]
\begin{center}
   \epsfig{file=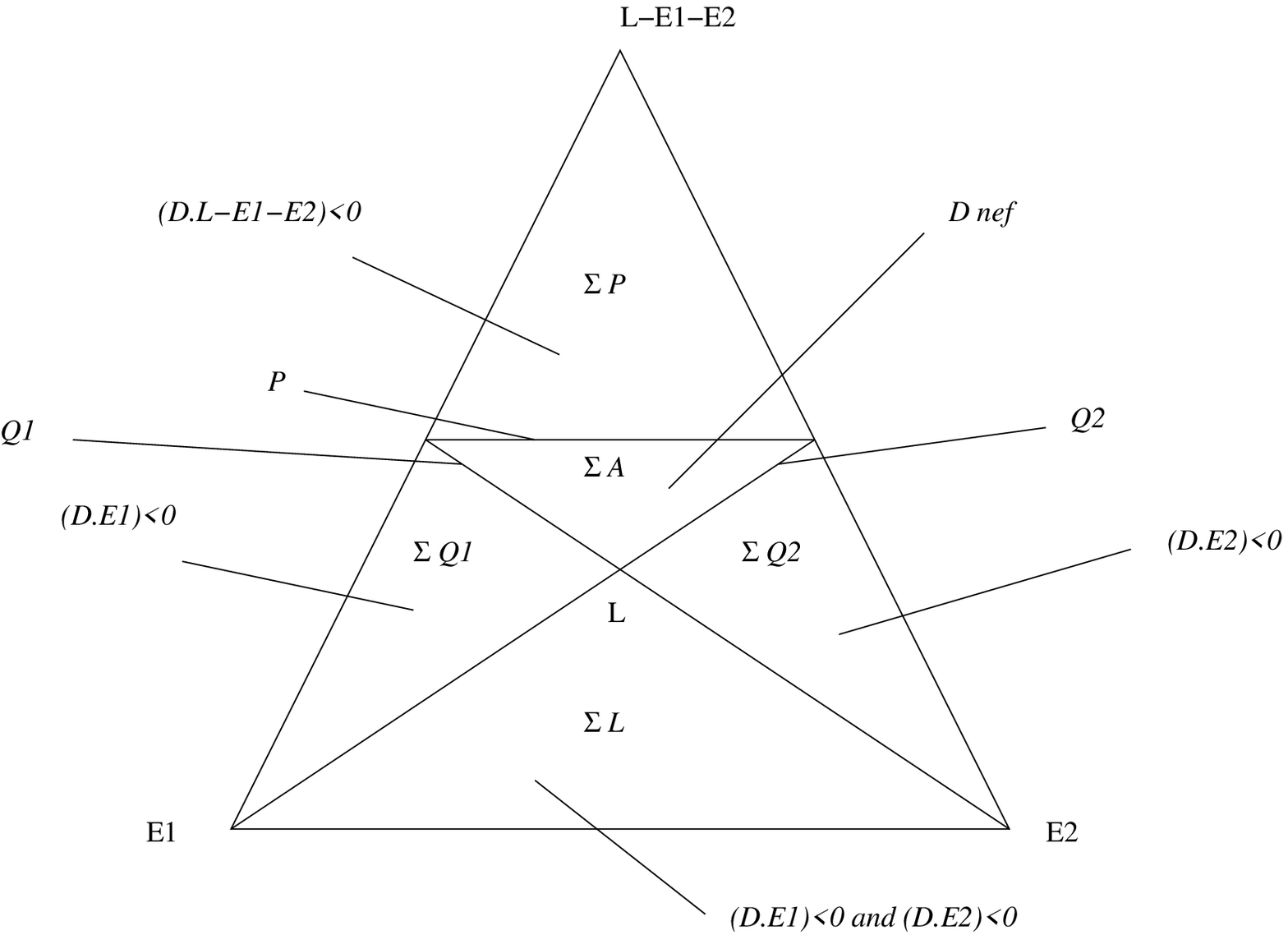,height=4in,width=4in} \label{fig1}
\end{center}

    Observe that not all possible combinations of negative divisors  occur. 
    This in part is accounted for by the fact that certain faces of the 
    nef cone are not contained in the big cone.
\end{example}

   Next we move on to spell out the connection with the Weyl action.
   To each root $\alpha$ one can associate a reflection of the lattice Pic($X$):
   $$
      \sigma_{\alpha}(D)=D+(D\cdot \alpha)\alpha.
   $$
   As a linear
   automorphism of Pic($X$), every $\sigma_{\alpha}$ descends uniquely
to the N\'eron-Severi space. The  group generated by the
reflections $\sigma_{\alpha}$, $\alpha$ a root, is called the {\em
Weyl group} $W(X)$ of the surface $X$. However, there is a much
smaller set of generators. For $r\geq 3$ the  roots
   $$
      \alpha_1 = L-E_1-E_2-E_3,\ \alpha_2 = E_2-E_1,\ \dots ,\ \alpha_r=E_r-E_{r-1}
   $$
   called {\em simple } roots already generate $W(X)$.
   The sets of big classes which intersect the same set of roots positively are called {\em Weyl chambers}.

\begin{proposition}
  Given a del Pezzo surface $X$, $W(X)$ is the set of automorphisms 
  of $\Pic(X)$ that leave $K_X$ fixed. It is finite for $r\leq 8$ 
  and acts transitively on $\II$ (for $r\geq 3$) and on ${\mathcal R}$ (for $r\geq 2$).
\end{proposition}

For a proof the reader is referred to \cite{demz}. The effect of the Weyl group on the volume and Zariski decompositions is given in the following

\begin{proposition}
   Let $D,D'$ be two big divisors that belong to the same Zariski chamber,
   and let $\sigma\in W(X)$ be an element of the Weyl group. Then:
   \begin{items}
      \item[(i)] The supports of the negative parts of $\sigma (D)$ and $\sigma (D' )$
      are also equal, ie. $\sigma (D)$ and $\sigma (D')$ also belong to the same Zariski chamber.
      \item[(ii)] $\vl{\sigma(D)}=\vl{D}$.
   \end{items}
\end{proposition}

\proof
   Assume $r\geq 3$, as the remaining cases are  easy to check. 
   Then $\sigma$ permutes the exceptional elements hence it takes 
   nef divisors to nef divisors. It is enough to check this statement 
   on the generating reflections so we can assume that ${\sigma}^2=1$. 
   If $C,D$ are  any divisors then
   \[
       \zj{ \sigma(D)\cdot C}=\zj{ D\cdot\sigma(C)}
   \]
   as $\sigma$ preserves the intersection form. By the previous proposition
   \[
       \sigma(\II) = \II
   \]
   so $D$ is nef if and only if $\sigma (D)$ is.

   Consequently, forming the \zd of $D$ commutes with the action of $\sigma$,
   i.e. if $D=P+\sum_{i=1}^{s}{a_iN_i}$ is the \zd of $D$
   then
   $$
      \sigma (D) = \sigma (P)  + \sum_{i=1}^{s}{a_i\sigma (N_i)}
   $$
   is the \zd of $\sigma (D)$. This proves both statements.
\endproof

\begin{remark}\rm
   We observe that for $r\geq 3$  the Weyl chambers and the volume chambers coincide.
\end{remark}

\subsection{K3 surfaces}

   The interplay between the volume and the Weyl action on del Pezzo surfaces is
   in some sense not typical. As we will see, on $K3$ surfaces,
   the volume function in not invariant under the action of the Weyl group.
   However, it is  still true that inside the big cone the volume chambers coincide with the Weyl chambers.

Let $Y$ be a projective $K3$ surface. Then there are no
$(-1)$-curves on $Y$ and the Mori cone is generated by either the nef cone itself  or by the rational curves with self-intersection $0$ and $-2$ (see \cite{kov}). In the first case there is nothing to prove, as every big divisor is ample. In the second case, a big divisor is nef if and only if it intersects every $(-2)$-curve
non-negatively. Hence we can restrict our attention to $(-2)$-curves.

To every $(-2)$-curve $E$ one associates as before the reflection
   $$
      {\sigma}_E(D)=D+(D\cdot E)E
   $$
of the N\'eron--Severi space. These elementary reflections then generate the Weyl group $W(Y)$ of the surface. The Weyl chamber consisting of divisor classes intersecting every $(-2)$-curve nonnegatively is called the fundamental Weyl chamber. In our case, it coincides with the nef cone. For terminology and basic results on Weyl groups of $K3$ surfaces, root systems and related topics the reader should consult \cite{dolg,root}.

   First, we show that contrary to del Pezzo surfaces $\vl{D}\not =\vl{\sigma (D)}$ in general.
   Let $P$ be  nef divisor on $Y$, $E$ a $(-2)$-curve such that $(P\cdot E)\not =0$. The \zd of ${\sigma}_E(P)$ is
\begin{equation}
 {\sigma}_E(P) = \zj{ P+\frac{1}{2}(P\cdot E)E} + \frac{1}{2}(P\cdot E)E\ .
\end{equation}
Hence for the volume of ${\sigma}_E(P)$ we obtain
\begin{eqnarray*}
 \vl{{\sigma}_E(P)}  & = & \vl{P+\frac{1}{2}(P\cdot E)E} \\
 & = & \si{P+\frac{1}{2}(P\cdot E)E}{2} \\
 & = & P^2 + {(P\cdot E)}^2 -\frac{1}{2}(P\cdot E) \\
 & = & \vl{P} +{(P\cdot E)}^2-\frac{1}{2}(P\cdot E)\ .
\end{eqnarray*}
 As $(P\cdot E)\not = 0$, this is not equal to $\vl{P}$ (one can replace $P$ by a sufficiently high multiple if necessary).

\begin{proposition}
For any $K3$ surface $Y$, the volume chambers and the Weyl chambers in the big cone are the same.
\end{proposition}

\proof
Let $D$ be a big divisor on $Y$. By definition $D$ is nef if and only if it belongs to the fundamental Weyl chamber. Therefore we will assume that $D$ is not nef.

Observe  that
\[
  {\Phi}_D\stackrel{\rm def}{=} \left\{ E \ | \textrm{ $E$ is a $(-2)$-curve with  $D\cdot E<0$ } \right\}
\]
   is finite. Also, ${\Phi}_D \subseteq \Neg(D)$  as $D.E<0$ for
   all $E\in {\Phi}_D$, hence in particular  $\left| {\Phi}_D
\right|=r \leq  \rho =\textrm{rank }NS(Y)$.
Again, by \cite{looi}, 1.11  we observe that $D=w(P)$, with $P$ nef, $w\in W(Y)$ and  $w={\sigma}_{E_1}\dots {\sigma}_{E_r}$  a minimal decomposition of $w$ in terms of the reflections corresponding to $(-2)$-curves and ${\Phi}_D=\left\{ E_1, \dots , E_r\right\}$. But then  ${\Phi}_D\supseteq \Neg(D)$ which implies  ${\Phi}_D=\Neg(D)$.  As Weyl chambers are characterized by ${\Phi}_D$ and volume chambers inside the big cone are characterized by $\Neg(D)$, we can conclude that indeed, for big divisors,  every volume chamber is a Weyl chamber and vice versa. \\
\endproof

\subsection{An example where the volume is not locally polynomial}\label{volume example}

   In this section we construct an example of
   a smooth projective threefold whose associated
   volume function is not given  locally by a polynomial.
   We use Cutkosky construction as recalled in Example \ref{nonrational}.
   We keep the notation introduced in that example.

   We may choose
   $$
      D=f_1+f_2\ , \ H=3(f_2+\delta)
   $$
   two ample divisors with classes $d$ and $h$ in such a way
   that the ray $d -\lambda h$ in $N^1(\RR)$ intersects the
   boundary of the nef cone at a point $\lambda=\sigma$.
   Take $A_1=D$ and $A_2= -H$, and define
   $\cale=\calo_S(A_1)\oplus \calo_S(A_2)$
   and $X=\bbP(\cale)$. Then
   $$
      \hh{0}{X}{{\OO}_X(k)} = \sum_{i+j=k}{\hh{0}{S}{{\OO}_S(iD-jH)}}
   $$
   by Lemma \ref{ruled threefold} and
   $$
      \hh{0}{S}{{\OO}_S(iD-jH)} =  \frac{1}{2} \si{id -jh}{2}
   $$
   if $\frac{j}{i}<\sigma $  and  $0$ if  $ \frac{j}{i} >\sigma $
   by Riemann-Roch on the abelian surface $S$.

   It is important to start with a non-ample line bundle $L={\OO}_{{\PP ({\mathcal E})}}(1)$ on $X$
   as inside the ample cone the volume is
   a polynomial function given by the self-intersection.
   Next, we will perturb $L$ by a small $\QQ$-divisor depending on a parameter
   $\epsilon$ and establish that the dependence of the volume on $\epsilon $ is not polynomial.
   Take
\[
   A_1(\ep)=D+\ep f_1\ ,\ A_2(\ep)=-H+\ep f_1 \ ,
\]
then
\[
   X={\PP}\left( {\OO}_S(A_1)\oplus {\OO}_S(A_2)\right)\simeq
   {\PP}\left( {\OO}_S(A_1(\ep)\oplus {\OO}_S(A_2(\ep))\right)
\]
for all $\ep$ and
\[
   L(\ep)\stackrel{\rm def}{=}{\OO}_{{\PP}\left( {\OO}_S(A_1(\ep)\oplus {\OO}_S(A_2(\ep))\right)}(1)
   \simeq L\otimes {\pi}^*{\OO}_S(\ep f_1)\ .
\]
Then we have
\[
   \hh{0}{S}{{\OO}_S(iD-jH+(i+j)\ep f_1)}= \frac{1}{2}\si{i\delta -jh+(i+j)\ep f_1}{2}
\]
if $ \frac{j}{i}< {\sigma}(\ep)$ and $ 0$ otherwise,  where ${\sigma}(\ep)= \frac{9+5\ep-\sr}{18-12\ep}$.
Put
\[
   q(x)=\si{xd-(1-x)h+\ep f_1}{2}\ ,
\]
   then as in \cite{PAG} we have
\[
 \vl{L(\ep)}=3\int_{\frac{1}{1+{\sigma}(\ep)}}^{1}{q(x)dx}\ .
\]

By plugging in the data of our example we obtain
\begin{eqnarray*}
\vl{L(\ep)} & = & \frac{1}{ {\left( - 27+7\ep+\sr \right) }^3 } \times \\
 && \times \left(  33480\ep+43128{\ep}^2+8748- 1692\sr + 14120{\ep}^3 \right. \\
 && \left. -3300\ep\sr -2740{\ep}^2\sr  \right. \\
 && \left. + 84{\ep}^3\sr +588{\ep}^4 \right)
\end{eqnarray*}
which is not a polynomial function of $\ep$.

   By taking a product of $X$ with a projective space of appropriate dimension 
   we can obtain examples in every dimension at least three: for $n\geq 3$ take 
   $Y=X\times {\PP}^{n-3}$, $M_{\ep}={\pi}_1^*L_{\ep}\otimes {\pi}_2^*H$, where 
   $H$ is the hyperplane class in ${\PP}^{n}$. Then by the K\"unneth formula for 
   the volume
   \[
      \vol{Y}{M_{\ep}}={n \choose 3}\vol{X}{L_{\ep}}\cdot \vol{{\PP}^{n-3}}{H} = 
      {n \choose 3}\vol{X}{L_{\ep}}\ .
   \]
   Therefore the volume function associated to $Y$ is not locally polynomial.

   The computations in
   this section were done with the help of the computer algebra
   package {\em Maple}.


\section{Auxiliary results}

\begin{lemma}\label{lemma matrix}
   Let $S$ be a negative definite $r\times r$-matrix over the reals such that
   $s_{ij}\ge 0$ for all $i\ne j$. Then all entries of the inverse
   matrix $\Sinv$ are $\le 0$.
\end{lemma}

   While we feel that the statement must be well-known, we
   indicate
   a proof for lack of a reference.

\proof
   We argue by induction on $r$. The case $r=1$ being clear, we
   assume $r\ge 2$. We have
   $$
      \Sinv=\frac1{\det S}S\adj=\frac1{(-1)^r\cdot|\det S|}S\adj
   $$
   so we need to show that all entries of the cofactor matrix
   $S\adj$ have the sign $(-1)^{r-1}$.
   Denoting by $S_{ij}$ the matrix obtained from $S$ by deleting
   the $i$-th row and the $j$-th column, the assertion is that
   the numbers $(-1)^{i+j}\det S_{ij}$ have the sign $(-1)^{r-1}$.
   For $i=j$ this is certainly true, since $S_{ii}$ is negative definite.
   By symmetry it is then sufficient to consider the case $i<j$.
   Expanding the determinant $\det S_{ij}$ with respect to the
   $i$-th column, we get an expression
   $$
      \det S_{ij}=\sum_{k=1}^{r-1}(-1)^{i+k}\det((S_{ij})_{ki})\cdot
      c_k  \eqno(*)
   $$
   where $c_k$ is a non-negative number.
   The essential point is now that one has
   $$
      \det((S_{ij})_{ki})=\det((S_{ii})_{k,j-1}) \ .
   $$
   The claim follows upon
   using this relation in equation $(*)$ above and applying the
   induction hypothesis.
\endproof

\begin{lemma}\label{curves in negative}
   Let $N$ be a negative definite divisor with irreducible
   components $E_1,\dots E_m$. Then the only irreducible curves
   in the linear span $\vecspan{E_1,\dots,E_m}$ in $\NR(X)$ are
   $E_1,\dots,E_m$.
\end{lemma}

\proof
   Let $C$ be any irreducible curve in
   $\vecspan{E_1,\dots,E_m}$.
   After re-indexing we can write
   $$
      C=\sum_{i=1}^k\alpha_iE_i-\sum_{i=k+1}^m\beta_iE_i
   $$
   where $\alpha_i>0$ and $\beta_i\ge0$.
   We must have $k\ge 1$, hence
   $$
      0>\(\sum_{i=1}^k\alpha_iE_i\)^2
      =\(\sum_{i=1}^k\alpha_iE_i\)\cdot\(C+\sum_{i=k+1}^m\beta_iE_i\) \ ,
   $$
   and we conclude that $C$ is one of the curves $E_i$.
\endproof

\begin{lemma}\label{missing lemma}
   Let $P$ be a nef and big divisor. Then any non-zero combination of
   curves in $\Null(P)$ is a negative divisor i.e. a divisor with
   negative definite intersection form.
\end{lemma}

\proof
   As $P^2>0$, this follows directly from the Hodge Index Theorem.
\endproof
   Now we recall, following closely Section 2.3.B in
   \cite{PAG}, the following construction of a ruled threefold
   underlying the construction of Cutkosky \cite{cut}. We need
   this for our Example \ref{nonrational} and in Section
   \ref{volume example}.

   Let $S$ be an irreducible projective surface, $A_1$, $A_2$
integral Cartier divisors on it. Put
\[ {\mathcal E}={\OO}_S(A_1)\oplus {\OO}_S(A_2)\ , \] and let \[ X={\PP}({\mathcal E})\ {\rm and\ } L={\OO}_{{\PP}({\mathcal E})}(1)\ .\]
   There is a close relation between the the properties of
   $L={\OO}_X(1)$ and those of $A_1$ and $A_2$.
\begin{lemma}\label{ruled threefold}
   With notation as above,
   \begin{items}
      \item[(i)] One has $\HH{0}{X}{L^{\otimes k}}={\bigoplus}_{a_1+a_2=k}{\HH{0}{S}{{\OO}_S(a_1A_2+a_2A_2)}}.$
      \item[(ii)] $L$ is ample if and only if $A_1$ and $A_2$ are ample on $S$.
      \item[(iii)] $L$ is nef if and only if $A_1$ and $A_2$ are nef on $S$.
      \item[(iv)] $L$ is big if and only if some $\NN$-linear combination of the $A_i$'s is a big divisor on $S$.
      \item[(v)] For $m\in \NN$, $L^{\otimes m}$ if free if and only if $mA_1$ and $mA_2$ are both free on $S$.
   \end{items}
\end{lemma}

In connection with the previous example, we will need the following interesting fact about the volume function not 
established in the literature so far.

\begin{proposition}[K\"unneth formula for the volume] Let $X_1,X_2$ be irreducible 
projective varieties of dimensions $n_1,n_2$, $L_1,L_2$ line bundles on the respective spaces. Then
\[
\vol{X_1\times X_2}{{\pi}_1^*L_1\otimes {\pi}_2^*L_2} = {n_1+n_2 \choose n_1} \vol{X_1}{L_1}\cdot
\vol{X_2}{L_2}\ .
\]
\end{proposition}
\proof
The statement follows from the K\"unneth formula for sheaves and the fact that 
\[
\vol{X}{L} = \lim_m \frac{\hh{0}{X}{mL}}{m^n/n!}
\]
in general.
\endproof



\enlargethispage{1cm}
\bigskip
   Tho\-mas Bau\-er,
   Fach\-be\-reich Ma\-the\-ma\-tik und In\-for\-ma\-tik,
   Philipps-Uni\-ver\-si\-t\"at Mar\-burg,
   Hans-Meer\-wein-Stra{\ss}e,
   D-35032~Mar\-burg, Germany.

\nopagebreak
   E-mail: {\tt tbauer@mathematik.uni-marburg.de}

\nopagebreak
\medskip
   Alex K\"uronya,
   Department of Mathematics,
   University of Michigan,
   Ann Arbor,
   MI~48109, USA

\nopagebreak
   E-Mail: {\tt akuronya@umich.edu}

\nopagebreak
\medskip
   Tomasz Szemberg,
   Akademia Pedagogiczna,
   Instytut Matematyki,
   ul. Podchor\c a\.zych 2,
   PL-30-084~Krak\'ow;
   Poland
   and
   Universit\"at Duisburg-Essen,
   Fachbereich 6 Mathematik,
   D-45117 Essen,
   Germany

\nopagebreak
   E-mail: {\tt szemberg@ap.krakow.pl}


\end{document}